\documentclass[12pt, twoside, eqno]{article}
\usepackage{latexsym}
\usepackage{amssymb}
\usepackage{amsfonts}
\begin{document}\bigskip

\noindent{\Large\bf From semigroups to $\Gamma$-semigroups and to 
hypersemigroups}\bigskip

\noindent{\bf Niovi Kehayopulu}\bigskip

\noindent August 10, 2016\bigskip{\small

\noindent{\bf Abstract.} This paper serves as an example to show the 
way we pass from semigroups to $\Gamma$-semigroups and to 
hypersemigroups.\smallskip

\noindent 2010 AMS Subject Classification: 20M99\\
Keywords: $\Gamma$-semigroup; hypersemigroup; regular; right (left) 
simple; bi-ideal}
\section{Introduction} The present paper is based on the paper  
``Note on bi-ideals in ordered semigroups" by N. Kehayopulu, J. S. 
Ponizovskii and M. Tsingelis in [2], and its aim is to show how we 
pass from semigroups to $\Gamma$-semigroups and to hypersemigroups. 
Many results on semigroups are transferred to $\Gamma$-semigroups 
just putting a ``Gamma" in the appropriate place. Many results on 
semigroups are transferred to hypersemigroups just replacing the 
multiplication ``$.$" of the semigroup by the operation ``$*"$ of the 
hypersemigroup. So, for any paper on $\Gamma$-semigroups or on 
hypersemigroups we have to indicate inside the paper the 1--2 papers 
on semigroups on which our paper is based. Even if these papers are 
cited in the References (many times are not even cited), they are 
always cited among a large amount 20--36 other published papers and, 
clearly, this is not enough. We say ``published papers" and this is 
because if we prove a result for a semigroup, we normally publish it, 
and do not keep its proof to transfer it to $\Gamma$-semigroups or to 
hypersemigroups or both. That's why the papers on semigroups (or 
ordered semigroups) must have been properly cited in any paper on 
$\Gamma$-semigroups or on hypersemigroups. In addition, there is 
another interesting information about the hypersemigroups, we will 
deal with in another paper.

For an ordered semigroup $(S,.,\le)$ and a subset $A$ of $S$, we 
denote by $(A]$ the subset of $S$ defined by: $(A]:=\{t\in S \mid 
t\le h \mbox { for some } a\in H\}$. A nonempty subset $A$ of $S$ is 
called a left (resp. right) ideal of $S$ if (1) $SA\subseteq A$ and 
(2) if $a\in A$ and $S\ni b\le a$, then $b\in A$, that is, if 
$(A]=A$. A nonempty subset $A$ of an ordered semigroup $S$ is called 
a bi-ideal of $S$ if (1) $ASA\subseteq A$ and (2) if $a\in A$ and 
$S\ni b\le a$, then $b\in A$. An ordered semigroup $S$ is called left 
(resp. right) simple if $S$ is the only left (resp. right) ideal of 
$S$. That is, if $A$ is a left (resp. right) ideal of $S$, then 
$A=S$. An ordered semigroup $S$ is called regular if for every $a\in 
S$ there exists $x\in S$ such that $a\le axa$. Equivalently, if $a\in 
(aSa]$ for every $a\in S$ or if $A\subseteq (ASA]$ for any subset $A$ 
of $S$. We have seen in [2], that an ordered semigroup $S$ is left 
(resp. right) simple if and only if, for every $a\in S$, we have 
$(Sa]=S$ (resp. $(aS]=S$). We have seen, that if an ordered semigroup 
is left and right simple, then it is regular and that in a regular 
ordered semigroup $S$ the bi-ideals and the bi-ideals which are at 
the same time subsemigroups of $S$ (i.e. the subidempotent bi-ideals 
of $S$) are the same. We have also proved that an ordered semigroup 
is both left simple and right simple if and only if does not contain 
proper bi-ideals. The result given in [2] can be transferred to 
$\Gamma$-semigroups just putting a ``Gamma" in the appropriate place. 
As far as the hypersemigroups is concerned, except of the Theorem 
3.12, where some little change is needed (as both the operation and 
the hyperoperation play a role in it), the rest are transferred from 
semigroups just replacing the multiplication of the semigroup in [2] 
by the operation ``$*$" of the hypersemigroup. For convenience, let 
us give a complete details of our arguments to justify what we say.
\section{On bi-ideals on $\Gamma$-semigroups} For two nonempty sets 
$M$ and $\Gamma$, we denote by $A\Gamma B$ the set containing the 
elements of the form $a\gamma b$ where $a\in A$, $\gamma\in\Gamma$ 
and $b\in B$. That is, we define $$A\Gamma B:=\{a\gamma b \mid a\in 
A, b\in B, \gamma\in\Gamma\}.$$Then $M$ is called a {\it 
$\Gamma$-semigroup} [1] if the following assertions are satisfied: 
\begin{enumerate}
\item[(1)] $M\Gamma M\subseteq M$.
\item[(2)] $a\gamma (b\mu c)=(a\gamma b)\mu c$ for all $a,b,c\in 
M$ and all $\gamma,\mu\in\Gamma$. 
\item[(3)] If $a,b,c,d\in M$ and $\gamma,\mu\in\Gamma$ such that 
$a=c$, $\gamma=\mu$ and $b=d$, then $a\gamma b=c\mu 
d$.\end{enumerate}
For $A=\{a\}$ we write, for short, $a\Gamma B$ instead of 
$\{a\}\Gamma B$ and for $B=\{b\}$ we write $A\Gamma b$ instead of 
$A\Gamma \{b\}$.

Let $M$ be a $\Gamma$-semigroup. A nonempty subset $A$ of $M$ is 
called a {\it left} (resp. {\it right}) {\it ideal} of $M$ if 
$M\Gamma A\subseteq A$ (resp. $A\Gamma M\subseteq A)$. $A$ is called 
an {\it ideal} of $M$ if it is both a left and a right ideal of $M$. 
It is called a {\it subsemigroup} of $M$ if $A\Gamma A\subseteq A$. 
Clearly, if $A$ is an ideal of $M$, then it is a subsemigroup of $M$. 
A nonempty subset $A$ of $M$ is called a {\it bi-ideal} of $M$ if 
$A\Gamma M\Gamma A\subseteq A$. A bi-ideal $A$ of $M$ is called {\it 
subidempotent} if it is a subsemigroup of $M$. A left ideal, right 
ideal or a bi-ideal $A$ of $M$ is called {\it proper} if $A\not=M$. A 
$\Gamma$-semigroup $M$ is called {\it left} (resp. {\it right}) {\it 
simple} if $M$ does not contain proper left (resp. right) ideals, 
that is if $A$ is a left (resp. right) ideal of $M$, then $A=M$. For 
a subset $A$ of $M$, we denote by $(A]$ the subset of $M$ defined by 
$(A]=\{t\in M \mid t\le a \mbox { for some } a\in A\}$. $M$ is called 
{\it regular} if, for every $a\in M$ there exist $x\in M$ and 
$\gamma,\mu\in\Gamma$ such that $a\le a\gamma x\mu a$. \smallskip

\noindent{\bf Proposition 2.1.} {\it let H be a $\Gamma$-semigroup. 
The following are equivalent:

$(1)$ H is regular.

$(2)$ $a\in a\Gamma M\Gamma a$ for every $a\in M$.

$(3)$ $A\subseteq A\Gamma M\Gamma A$ for every $A\subseteq M$. } 
\medskip

\noindent{\bf Proof.} $(1)\Longrightarrow (2)$. Let $a\in M$. Since 
$H$ is regular, there exist $x\in M$ and $\gamma,\mu\in\Gamma$ such 
that $a\in a\gamma x\mu a\in a\Gamma M\Gamma a$.\\$(2)\Longrightarrow 
(3)$. Let $A\subseteq M$ and $a\in A$. By (2), we have $a\in a\Gamma 
M\Gamma a\subseteq A\Gamma M\Gamma A$.\\$(3)\Longrightarrow (1)$. Let 
$A\subseteq M$ and $a\in A$. Since $\{a\}\subseteq M$, by (3), we 
have $\{a\}\subseteq \{a\}\Gamma M\Gamma \{a\}$. Then $a\in a\Gamma 
M\Gamma a$, and there exist $x\in M$ and $\gamma,\mu\in\Gamma$ such 
that $a=a\gamma x\mu a$, so (1) holds. $\hfill\Box$\medskip

\noindent{\bf Proposition 2.2.} {\it Let M be a $\Gamma$-groupoid. If 
$M\Gamma a=M$ for every $a\in M$, then $M$ is left simple. 
``Conversely", if M is a $\Gamma$-semigroup and M is left simple 
then, for every $a\in M$, we have $M\Gamma a=M$}.\medskip

\noindent{\bf Proof.} $\Longrightarrow$. Let $T$ be a left ideal of 
$M$. Then $T=M$. Indeed: Let $a\in M$. Take an element $b\in T$ 
$(T\not=\emptyset)$. By hypothesis, we have $M=M\Gamma b\subseteq 
M\Gamma T\subseteq T$, so $M=T$.\\$\Longleftarrow$. Let $M$ be a left 
simple $\Gamma$-semigroup and $a\in M$. The set $M\Gamma a$ is a left 
ideal of $M$. Indeed, it is a nonempty subset of $M$ (as 
$M,\Gamma\not=\emptyset)$ and $M\Gamma (M\Gamma a)=(M\Gamma M)\Gamma 
a\subseteq M\Gamma a$. Since $M$ is left simple, we have $M\Gamma 
a=M$. $\hfill\Box$\medskip

\noindent{\bf Corollary 2.3.} {\it A $\Gamma$-semigroup M is left 
simple (resp. right simple) if and only if, for every $a\in M$, we 
have $M\Gamma a=M \mbox { (resp. } a\Gamma M=M)$}.\medskip

\noindent{\bf Proposition 2.4.} {\it Let M be a $\Gamma$-semigroup. 
If M is left simple and right simple, then M is regular. }\medskip

\noindent{\bf Proof.} Let $a\in M$. Since $M\Gamma a=M$ and $a\Gamma 
M=M$, we have $a\in A\Gamma M=a\Gamma(M\Gamma a)$, so $a\in a\Gamma 
M\Gamma a$ and, by Proposition 2.1, $M$ is regular. 
$\hfill\Box$\medskip

\noindent{\bf Proposition 2.5.} {\it In a regular $\Gamma$-semigroup, 
the bi-ideals and the subidempotent bi-ideals are the same.} 
\medskip

\noindent{\bf Proof.} Let $B$ be a bi-ideal of $M$. Then $B\Gamma 
M\Gamma B\subseteq B$. Since $M$ is regular, by Proposition 2.1, we 
have $B\subseteq B\Gamma M\Gamma B$, so we get $B=B\Gamma M\Gamma B$. 
Then we have$$B\Gamma B=(B\Gamma M\Gamma B)\Gamma (B\Gamma M\Gamma 
B)=B\Gamma (M\Gamma B\Gamma B\Gamma M)\Gamma B.$$Since 
\begin{eqnarray*}M\Gamma B\Gamma B\Gamma M&=&M\Gamma (B\Gamma 
B)\Gamma M\subseteq M\Gamma (M\Gamma M)\Gamma M\\&\subseteq& M\Gamma 
(M\Gamma M)\subseteq M\Gamma M\subseteq M,\end{eqnarray*} we have 
$B\Gamma B\subseteq B\Gamma M\Gamma B=B$, so $B$ is a subsemigroup of 
$M$. $\hfill\Box$\medskip

\noindent{\bf Proposition 2.6.} {\it Let M be a $\Gamma$-semigroup. 
If A is a left (or right) ideal of M, then A is a bi-ideal of 
M}.\medskip

\noindent{\bf Proof.} Let $A$ be a left ideal of $M$. Then $A\Gamma 
M\Gamma A=A\Gamma (M\Gamma A)\subseteq A\Gamma A$. Since $A$ is a 
subsemigroup of $M$, we have $A\Gamma A\subseteq A$, then $A\Gamma 
M\Gamma A\subseteq A$, and $A$ is a bi-ideal of $M$. 
$\hfill\Box$\medskip

\noindent{\bf Proposition 2.7.} {\it Let M be a $\Gamma$-semigroup. 
Then, for any nonempty subsets $A,B,C$ of $M$, we have

$(1)$ $(A\cup B)\Gamma C=A\Gamma C\cup B\Gamma C$.

$(2)$ $C\Gamma (A\cup B)=C\Gamma A\cup C\Gamma B$.}\smallskip

For a $\Gamma$-semigroup $M$ and an element $b$ of $M$, we denote by 
$L(b)$ (resp. $R(b)$) the left (resp. right) ideal of $M$ generated 
by $b$. \medskip

\noindent{\bf Proposition 2.8.} {\it Let M be a $\Gamma$-semigroup 
and $b\in M$. Then we have the following:

$(1)$ $L(b)=\{b\}\cup M\Gamma b$.

$(2)$ $R(b)=\{b\}\cup b\Gamma M$.}\smallskip

\noindent{\bf Proof.} The set $\{b\}\cup M\Gamma b$ is a subset of 
$M$ containing $b$. Moreover, it is a left ideal of $M$. 
Indeed:\begin{eqnarray*}M\Gamma {\Big(}\{b\}\cup M\Gamma 
b{\Big)}&=&M\Gamma b\cup (M\Gamma M)\Gamma b \mbox { (by Proposition 
2.7)}\\&=&M\Gamma b \mbox { (since} M\Gamma M\subseteq 
M)\\&\subseteq&M\Gamma {\Big(}\{b\}\cup M\Gamma 
b{\Big)}.\end{eqnarray*} And if $T$ is a left ideal of $M$ such that 
$b\in T$, then $\{b\}\cup M\Gamma b\subseteq T\cup M\Gamma T=T$. So 
the set $\{b\}\cup M\Gamma b$ is the least left ideal of $M$ 
containing $b$, that is, $L(b)=\{b\}\cup M\Gamma b$. The proof of (2) 
is similar. $\hfill\Box$\smallskip

Exactly as in the Proposition 1 in [2], just putting a ``Gamma" where 
is needed, we can prove the next theorem. For the sake of 
completeness, we will give its proof.\smallskip

\noindent{\bf Theorem 2.9.} {\it A $\Gamma$-semigroup M is left 
simple and right simple if and only if does not contain proper 
bi-ideals}. \smallskip

\noindent{\bf Proof.} $\Longrightarrow$. Let $A$ be a bi-ideal of 
$M$. Then $A=M$. In fact: Let $a\in M$. Take an element $b\in A$ 
$(A\not=\emptyset)$. We consider the left ideal of $M$ generated by 
$b$, that is the set $L(b)=\{b\}\cup M\Gamma b$. Since $M$ is left 
simple, we have $L(b)=M$. Since $a\in M$, we have $a\in L(b)$. Then 
$a=b$ or $a\in M\Gamma b$. If $a=b$ then, since $b\in A$, we have 
$a\in A$. Let $a\in M\Gamma b$. Then $a=x\gamma b$ for some $x\in M$, 
$\gamma\in\Gamma$. We consider the right ideal of $M$ generated by 
$b$, that is the set $R(b)=b\cup b\Gamma M$. Since $M$ is right 
simple, we have $R(b)=M$. Since $x\in M$, we have $x\in R(b)$. Then 
we have $x=b$ or $x\in b\Gamma M$. If $x=b$, then $a=x\gamma 
b=b\gamma b\in A\Gamma A\subseteq A$ (by Propositions 2.4 and 2.5). 
Let $x\in b\Gamma M$. Then $x=b\mu y$ for some $\mu\in\Gamma$, $y\in 
M$. Then we have $a=x\gamma b=b\mu y\gamma b\in A\Gamma M\Gamma 
A\subseteq A$, so $A=M$.\\$\Longleftarrow$. Let $A$ be a left ideal 
of $M$. By Proposition 2.6, $A$ is a bi-ideal of $M$. By hypothesis, 
we have $A=M$, so $M$ is left simple. In a similar way we prove that 
$M$ is right simple. $\hfill\Box$

\section{On bi-ideals in hypersemigroups}
An {\it hypergroupoid} is a nonempty set $H$ with an hyperoperation 
$$\circ : H\times H \rightarrow {\cal P}^*(H) \mid (a,b) \rightarrow 
a\circ b$$on $H$ and an operation $$* : {\cal P}^*(H)\times {\cal 
P}^*(H) \rightarrow {\cal P}^*(H) \mid (A,B) \rightarrow A*B$$ on 
${\cal P}^*(H)$ (induced by the operation of $H$) such that 
$$A*B=\bigcup\limits_{(a,b) \in\,A\times B} {(a\circ b)}$$for every 
$A,B\in {\cal P}^*(H)$. As the operation ``$*$" depends on the 
hyperoperation ``$\circ$", an hypergroupoid can be also denoted by 
$(H,\circ)$ (instead of $(H,\circ,*)$).

By the definition of the operation ``$\circ$", if $H$ is an 
hypergroupoid, then $$H*H\subseteq H.$$
If $H$ is an hypergroupoid then, for every $x,y\in H$, we have 
$\{x\}*\{y\}=x\circ y.$ An hypergroupoid $H$ is called {\it 
hypersemigroup} if $$\{x\}*(y\circ z)=(x\circ y)*\{z\}$$for every 
$x,y,z\in H$. So an hypergroupoid $H$ is an hypersemigroup if and 
only if, for every $x,y,z\in H$, we have 
$\{x\}*{\Big(}\{y\}*\{z\}{\Big)}={\Big(}\{x\}*\{y\}{\Big)}*\{z\}.$

If $H$ is an hypersemigroup then, for any nonempty subsets $A,B,C$ of 
$H$, we have $(A*B)*C=A*(B*C)=\bigcup\limits_{(a,b,c)\in A\times 
B\times C} {{\Big(}\{a\}*\{b\}*\{ c\}\Big) }$. Thus we can write 
$(A*B)*C=A*(B*C)=A*B*C$.

Since $({\cal P}^*(H),*)$ is a semigroup, for any  product $A_1*A_2* 
..... * A_n$ of elements of ${\cal P}^*(H)$ we can put parentheses in 
any place beginning with some $A_i$ and ending in some $A_j$ $(1\le 
i,j\le n)$.

The following proposition, though clear, plays an essential role in 
the theory of hypergroupoids.\medskip

\noindent{\bf Proposition 3.1.} {\it Let $(H,\circ)$ be an 
hypergroupoid, $x\in H$ and $A,B\in {\cal P}^*(H)$. Then we have the 
following:

$(1)$ $x\in A*B$ $\Longleftrightarrow$ $x\in a\circ b$ for some $a\in 
A$, $b\in B$.

$(2)$ If $a\in A$ and $b\in B$, then $a\circ b\subseteq 
A*B$.}\medskip

\noindent{\bf Proposition 3.2.} {\it Let H be an hypergroupoid and 
$A,B,C\in {\cal P}^*(H)$. Then $A\subseteq B$, implies $A*C\subseteq 
B*C$ and $C*A\subseteq C*B$. Equivalently,

$A\subseteq B$ and $C\subseteq D$ implies $A*C\subseteq B*D$ and 
$C*A\subseteq D*B$. } \medskip

\noindent{\bf Proof.} Let $A\subseteq B$ and $x\in A*C$. By 
Proposition 3.1, $x\in a\circ c$ for some $a\in A$ and $c\in C$. Then 
$x\in a\circ c$, where $a\in B$ and $c\in C$ thus, by Proposition 
3.1, we get $x\in B*C$. $\hfill\Box$\medskip

Let $H$ be an hypergroupoid. A nonempty subset $A$ of $H$ is called a 
{\it left} (resp. {\it right}) {\it ideal} of $H$ if $H*A\subseteq A$ 
(resp. $A*H\subseteq A)$. It is called a {\it subsemigroup} of $H$ if 
$A*A\subseteq A$. Clearly, the left ideals and the right ideals of 
$H$ are subsemigroups of $H$. A left (resp. right) ideal $A$ of $H$ 
is called {\it proper} if $A\not=H$. An hypergroupoid $H$ is called 
{\it left simple} (resp. {\it right simple}) if $H$ does not contain 
proper left (resp. right) ideals. That is, for any left (resp. right) 
ideal $A$ of $H$ we have $A=H$. Let now $H$ be an hypersemigroup. A 
nonempty subset $B$ of $H$ is called a {\it bi-ideal} of $H$ if 
$B*H*B\subseteq B$. A bi-ideal of $H$ is called {\it proper} if 
$B\not=H$. A subset $B$ of an hypergroupoid $H$ is called {\it 
subidempotent} if $B*B\subseteq B$, that is if it is a subgroupoid of 
$H$.\medskip

\noindent{\bf Definition 3.3.} An hypergroupoid $H$ is called {\it 
regular} if, for every $a\in H$, there exist $x\in H$ such that $a\in 
(a\circ x)*\{a\}$ or $a\in \{a\}*(x\circ a)$.

An hypersemigroup $H$ is called {\it regular} if, for every $a\in H$, 
there exist $x\in H$ such that $a\in (a\circ x)*\{a\}$. (If $H$ is an 
hypersemigroup then, clearly, $(a\circ x)*\{a\}=\{a\}*(x\circ 
a)=\{a\}*\{x\}*\{a\}$).\medskip

\noindent{\bf Proposition 3.4.} {\it Let H be an hypersemigroup. Then 
the following are equivalent:

$(1)$ H is regular.

$(2)$ $a\in \{a\}*H*\{a\}$ for every $a\in H$.

$(3)$ $A\subseteq A*H*A$ for every $A\in {\cal P}^*(H)$. }\smallskip

\noindent{\bf Proof.} $(1)\Longrightarrow (2)$. Let $a\in H$. Since 
$H$ is regular, there exists $x\in H$ such that $a\in 
\{a\}*\{x\}*\{a\}\subseteq \{a\}*H*\{a\}$.\\$(2)\Longrightarrow (3)$. 
Let $A\in {\cal P}^*(H)$ and $a\in A$. Since $a\in H$, by (2), we 
have\\ $a\in \{a\}*H*\{a\}\subseteq A*H*A$, thus we get $a\in 
A*H*A$.\\ $(3)\Longrightarrow (1)$. Let $a\in H$. Since $\{a\}\in 
{\cal P}^*(H)$. by (3), we have$$\{a\}\subseteq 
\{a\}*H*\{a\}=\bigcup\limits_{(c,x,e) \in \{ a\}  \times H \times \{ 
a\} } {{\Big(}\{ c\} *\{ x\} *\{ e\}{\Big)}}.$$Then there exists 
$(c,x,e)\in \{a\}\times H\times \{a\}$ such that $a\in 
\{c\}*\{x\}*\{e\}$. Since $c=a$, $e=a$ and $x\in H$, we have $a\in 
\{a\}*\{x\}*\{a\}$, where $x\in H$, so $H$ is regular. $\hfill\Box$
\medskip

\noindent{\bf Proposition 3.5.} {\it Let H be an hypergroupoid. If 
$H*\{a\}=H$ for every $a\in H$, then H is left simple. ``Conversely", 
if H is an hypersemigroup and H is left simple then, for every $a\in 
H$, we have $H*\{a\}=H$.}\medskip

\noindent{\bf Proof.} Suppose $H*\{a\}=H$ for every $a\in H$ and let 
$T$ be a left ideal of $H$. Then $T=H$. Indeed: Let $a\in H$. Take an 
element $b\in T$ $(T\not=\emptyset)$. Then we have 
$H=H*\{b\}\subseteq H*T\subseteq T$, so $T=H$. Let now $H$ be a left 
simple hypersemigroup and $a\in H$. The set $H*\{a\}$ is a left ideal 
of $H$. Indeed: Since $H, \{a\}\in {\cal P}^*(H)$, by the definition 
of ``$*$", we have $H*\{a\}\in {\cal P}^*(H)$, so $H*\{a\}$ is a 
nonempty subset of $H$. Moreover, $H*(H*\{a\})=(H*H)*\{a\}\subseteq 
H*\{a\}$. Since $H$ is left simple, we have $H*\{a\}=H$. 
$\hfill\Box$\\
The right analogue of Proposition 3.5 also holds and we have the 
following\medskip

\noindent{\bf Corollary 3.6.} {\it An hypersemigroup H is left simple 
(resp. right simple) if and only if, for every $a\in H$, we 
have$$H*\{a\}=H \mbox { (resp.}\; \{a\}*H=H).$$}{\bf Proposition 
3.7.} {\it Let H be an hypersemigroup. If H is left simple and right 
simple, then H is regular}. \medskip

\noindent{\bf Proof.} Let $a\in H$. Since $H$ is left simple, by 
Corollary 3.6, we have $H*\{a\}=H$. Since $H$ is right simple, we 
have $\{a\}*H=H$. Since $a\in H$, we have$$a\in 
\{a\}*H=\{a\}*(H*\{a\})=\{a\}*H*\{a\}.$$By Proposition 3.4, $H$ is 
regular. $\hfill\Box$\medskip

\noindent{\bf Proposition 3.8.} {\it In a regular hypersemigroup the 
bi-ideals and the subidempotent bi-ideals are the same.}\medskip

\noindent{\bf Proof.} Let $H$ be a regular hypersemigroup and $B$ a 
bi-ideal of $H$. Since $H$ is regular and $B\in {\cal P}^*(H)$, by 
Proposition 3.4, we have $B\subseteq B*H*B$. Since $B$ is a bi-ideal 
of $H$, we have $B*H*B\subseteq B$. Thus we have $B=B*H*B$. Then
$$B*B=(B*H*B)*(B*H*B)=B*(H*B*B)*H*B.$$Since $H*B*B\subseteq 
H*H*H\subseteq H$, we have $$B*B\subseteq B*(H*H)*B\subseteq 
B*H*B=B.$$ $\hfill\Box$\\
\noindent{\bf Proposition 3.9.} {\it Let H be an hypersemigroup. If A 
is a left (or right) ideal of H, then A is a bi-ideal of H.}\medskip

\noindent{\bf Proof.} Let $A$ be a left ideal of $H$. Then 
$H*A\subseteq A$. Then we have$$A*H*A=A*(H*A)\subseteq A*A\subseteq 
H*A\subseteq A,$$so $A$ is a bi-ideal of $H$. Similarly, the right 
ideals of $H$ are bi-ideals of $H$. $\hfill\Box$\medskip

\noindent{\bf Proposition 3.10.} {\it Let H be an hypergroupoid and 
$A,B,C\in {\cal P}^*(H)$. Then we have the following:

$(1)$ $(A\cup B)*C=(A*C)\cup (B*C)$.

$(2)$ $C*(A\cup B)=(C*A)\cup (C*B)$.}\medskip

\noindent{\bf Proof.} Let $x\in (A\cup B)*C$. Then $x\in a\circ b$ 
for some $a\in A\cup B$, $b\in C$. If $a\in A$, then $a\circ 
b\subseteq A*C$. If $a\in B$, then $a\circ b\subseteq B*C$. Thus we 
have $x\in (A*C)\cup (B*C)$. Let now $x\in (A*C)\cup (B*C)$. If $x\in 
A*C$, then $x\in a\circ c$ for some $a\in A$, $c\in C$, so $x\in 
(A\cup B)*C$. If $x\in B*C$, then $x\in b\circ c$ for some $b\in B$, 
$c\in C$, so $x\in (A\cup B)*C$. The proof of (2) is similar. 
$\hfill\Box$\medskip

For an hypergroupoid $H$ and an element $b$ of $H$, we denote by 
$L(b)$ (resp. $R(b)$) the left (resp. right) ideal of $H$ generated 
by $b$.\medskip

\noindent{\bf Proposition 3.11.} {\it Let H be an hypersemigroup and 
$b\in H$. Then we have the following:

$(1)$ $L(b)=\{b\}\cup {\Big(}H*\{b\}{\Big)}$.

$(2)$ $R(b)=\{b\}\cup {\Big(}\{b\}*H{\Big)}$.} \medskip

\noindent{\bf Proof.} (1) Let $b\in H$ and $T:=\{b\}\cup 
{\Big(}H*\{b\}{\Big)}$. Then $b\in T$ and $T\subseteq H\cup (H*H)=H$, 
so $T$ is a nonempty subset if $H$ containing $b$. Moreover, $T$ is a 
left ideal of $H$. Indeed:\begin{eqnarray*}H*T&=&H*{\Bigg(}\{b\}\cup 
{\Big(}H*\{b\}{\Big){\Bigg)}}\\&=&{\Big(}H*\{b\}{\Big)}\cup 
{\Bigg(}H*{\Big(}H*\{b\}{\Big)}{\Bigg)}\mbox { (by Proposition 
3.10)}\\&=&{\Big(}H*\{b\}{\Big)}\cup 
{\Big(}(H*H)*\{b\}{\Big)}\\&=&H*\{b\} \mbox { (since } H*H\subseteq 
H)\\&\subseteq&\{b\}\cup {\Big(}H*\{b\}{\Big)}=T. \end{eqnarray*}
Let now $K$ be a left ideal of $H$ such that $b\in K$. 
Then$$T=\{b\}\cup {\Big(}H*\{b\}{\Big)}\subseteq K\cup (H*K)=K,$$ so 
$L(b)=T$. The proof of (2) is similar. $\hfill\Box$\medskip

\noindent{\bf Theorem 3.12.} {\it An hypersemigroup H is both left 
simple and right simple if and only if H does not contain proper 
bi-ideals}.\medskip

\noindent{\bf Proof.} $\Longrightarrow$. Let $A$ be a bi-ideal of 
$H$. Then $A=H$. In fact: Let $a\in H$. Take an element $b\in A$ 
$(A\not=\emptyset)$. We consider the left ideal of $H$ generated by 
$b$, that is the set $L(b)=\{b\}\cup {\Big(}H*\{b\}{\Big)}$. Since 
$H$ is left simple, we have $L(b)=H$. Since $a\in H$, we have $a\in 
L(b)$. Then $a=b$ or $a\in H*\{b\}$. If $a=b$ then, since $b\in A$, 
we have $a\in A$. Let $a\in H*\{b\}$. Then, by Proposition 3.1, $a\in 
x\circ b$ for some $x\in H$. We consider the right ideal of $H$ 
generated by $b$, that is the set $R(b)=\{b\}\cup 
{\Big(}\{b\}*H{\Big)}$. Since $H$ is right simple, we have $R(b)=H$. 
Since $x\in H$, we have $x\in R(b)$. Then we have $x=b$ of $x\in 
\{b\}*H$. If $x=b$, then $a\in x\circ b=b\circ b\subseteq 
A*A\subseteq A$. Let $x\in \{b\}*H$. Then, by Proposition 3.1, we 
have $x\in b\circ y$ for some $y\in H$. By Proposition 3.2, we have 
$$a\in x\circ b\subseteq (b\circ y)*\{b\}=\{b\}*\{y\}*\{b\}\subseteq 
A*H*A\subseteq A,$$ so $a\in A$.\\$\Longleftarrow$. Let $A$ be a left 
ideal of $H$. Then $A$ is a bi-ideal of $H$. By hypothesis, we have 
$A=H$, so $H$ is left simple. In a similar way we prove that $H$ is 
right simple. $\hfill\Box$

An hypersemigroup $H$ is said to be {\it hypergroup} if the following 
assertions are satisfied:

(1) there exists $e\in H$ such that $a\circ e=e\circ a=\{a\}$ for 
every $a\in H$ and

(2) for every $a\in H$ there exists $a^{-1}\in H$ such that $a\circ 
a^{-1}=a^{-1}\circ a=\{e\}$.\medskip

\noindent{\bf Theorem 3.13.} {\it Let H be an hypergroup. Then H does 
not contain proper bi-ideals}.\medskip

\noindent{\bf Proof.} Let $A$ be a bi-ideal of $H$ and $a\in H$. Take 
an element $b\in A$ $(A\not=\emptyset)$. Since $b\in H$, there exists 
$b^{-1}\in H$ such that $b\circ b^{-1}=b^{-1}\circ b=\{e\}$, where 
$e\in H$ such that $x\circ e=e\circ x=\{x\}$ for every $x\in H$. We 
have\begin{eqnarray*}\{a\}&=&a\circ e\subseteq (e\circ 
a)*\{e\}=\{e\}*\{a\}*\{e\}\\&=&(b\circ b^{-1})*\{a\}*( b^{-1}\circ 
b)\\&=&\{b\}*{\Big(}\{b^{-1}\}*\{a\}*\{b^{-1}\}{\Big)}*\{b\}\\
&\subseteq&A*H*A\subseteq A.\end{eqnarray*}Thus we have $a\in A$, so 
$A=H$. $\hfill\Box$\medskip

\noindent{\bf Problem.} Find an example of an hypersemigroup which 
does not contain proper bi-ideals and it is not an hypergroup. 
{\small\bigskip

\medskip

\noindent University of Athens\\
Department of Mathematics\\
15784 Panepistimiopolis, Greece\\
email: nkehayop@math.uoa.gr

\end{document}